\documentclass[12pt]{amsart}
\usepackage{amsfonts,amssymb,eucal}
\begin{document}
                                        
\newcommand{\C}{{\mathbb C}}
\newcommand{\half}{{\textstyle{\frac{1}{2}}}}
\newcommand{\bL}{{\overline{L}}}
\newcommand{\cL}{{\mathcal L}}
\newcommand{\BL}{{\widetilde{L}}}
\newcommand{\bM}{{\overline{\mathcal M}}}
\newcommand{\Gal}{{\rm Gal}}
\newcommand{\PP}{{\mathbb P}}
\newcommand{\bpi}{{\mbox{\boldmath{$\pi$}}}}
\newcommand{\bPi}{{\mbox{\boldmath{$\Pi$}}}}
\newcommand{\rel}{{\; {\rm rel} \;}}                             
\newcommand{\Q}{{\mathbb Q}}
\newcommand{\R}{{\mathbb R}}
\newcommand{\Z}{{\mathbb Z}}

\title {Notes on the geography of the plane at infinity}\bigskip

\author{Jack Morava}

\address{Department of Mathematics, Johns Hopkins University,
Baltimore, Maryland 21218}

\email{jack@math.jhu.edu}

\thanks{The author was supported by the NSF}

\date {1 October 2006}

\begin{abstract}
This note is a kind of sequel to [19]: the questions behind it have their start in early work
of Kapranov on Drinfel'd's asymptotic zones for the KZ equations, and they grew under the
influence of Consani and Marcolli's insight that completely decomposed curves are natural
basepoints in Arakelov geometry. It has benefitted from conversations with J. Baez, S.
Devadoss, M. Davis, P. Etingof, A. Henriques, T. Januszkiewicz, Y. Manin, J. Stasheff, D.
Stevenson, and others, but it remains incomplete and inconclusive. I have compiled it as
evidence of a rich and mysterious structure, some kind of fundamental groupoid at infinity on the
projective line, whose existence is coming slowly into focus. \end{abstract}

\maketitle

\section{Variations on the fundamental group} \bigskip

\noindent                      
Given a nice space $A$ with a finite subset $B$ of base points, the classical fundamental
groupoid $\bpi \; (A \rel B)$ is the category with elements of $B$ as objects, and equivalence
classes of paths, e.g. from $v \in B$ to $v' \in B$, as morphisms between them.
\bigskip
         
\noindent  
{\bf 1.1} If $A \subset X$ is contained in some reasonable larger space, this construction can be
elaborated to define a {\bf two}-groupoid (ie, a two-category in which all morphisms and
two-morphisms are invertible). If for simplicity $B$ consists of a single basepoint, and $A$ is
connected, then $\bPi \; (X \rel A)$ will have one object, with $\pi_1(A,b)$ as its self-maps. The 
boundary homomorphism
\[
\partial : \pi_2(X,A,b) \to \pi_1(A,b)
\]
in the homotopy exact sequence of a pair defines a two-category with  
\[
\{ \phi \in \pi_2(X,A,b) \; | \; \partial (\phi) = \sigma_1^{-1} \cdot \sigma_0 \in \pi_1(A,b) \}
\]
as two-morphisms between $\sigma_0,\sigma_1 \in \pi_1(A,b)$. In this two-category,
paths in $A$ are the morphisms, and the two-morphisms are `bubbles' in $X$ bounding
these paths; to check the composition axioms it is probably easiest to think of $\pi_n(X,A,b)$ as
defined by maps of the $n$-cube, relative to its boundary, to $(X,A)$. This definition continues
to make sense when $B$ has more than one element. \bigskip

\noindent
{\bf 1.2} The case in which $A = V(\R)$ and $X = V(\C)$ are the real and complex points of an
algebraic variety seems particularly interesting: real bubbles (i.e. morphisms $\beta : \PP_1 \to
V$ defined over $\R$) have invariants in the second relative homotopy group. The Galois group
of $\C$ over $\R$ defines an involution on $V(\C)$ with $V(\R)$ as fixed points; restricting a
map
\[
(\PP_1(\C),\PP_1(\R)) \to (V(\C),V(\R)) 
\]
to the lower hemisphere of $\PP_1(\C)$ defines an element 
\[
\deg (\beta) \in \pi_0 \; {\rm Maps}(D^2, \partial D^2;V(\C),V(\R)) = \pi_2(V(\C),V(\R)) 
\]
based at $0 \in \PP_1$. This group fits in an exact sequence
\[
\cdots \to \pi_2 V(\C) \to \pi_2(V(\C),V(\R)) \to \pi_1V(\R) \to \cdots \;;
\]
complex conjugation $\beta \mapsto \beta^c$ in the domain acts trivially in the group on the
right, and by reversing signs on the left. The class $\deg(\beta - \beta^c)$ lifts to $\pi_2 V(\C)$,
where it represents the class of $\beta_*[\PP_1(\C)]$. \bigskip

\noindent
One of the purposes of this note is to propose
\[
\bPi \; (V(\C) \rel V(\R) \rel F) \; := \; Q\pi_1 (V \rel F)
\]
as an {\it ad hoc} kind of {\bf quantum} fundamental (two-)group(oid) of $V$ (with respect to
some suitable collection $F \subset V(\R)$ of basepoints). \bigskip

\noindent
Projective space is an interesting example: when $n>1, \; \bPi \; (\PP_n,\infty)$ has one object,
one morphism, and $\Z$ as its two-morphisms; but $\bPi \; (\PP_1,\infty)$ has one object, $\Z$
as its morphisms, and $2\Z$ as two-morphisms. In both cases the Galois group of $\C$ over $\R$
acts as $\pm 1$ on the two-morphisms, trivially on the morphisms. \bigskip

\noindent 
When $V$ is defined over $\Q$ or $\Z$, the set $V(\Q)$ presents itself as a potential supply of
basepoints for this construction; in the case of $\PP_1$, the points $0,1$ and $\infty$ stand
out. There is an algebraic version of the fundamental groupoid, which involves the specification
of suitable tangent vectors at the base points, cf. eg. [14], but I will neglect that refinement here. I
suspect that related higher invariants of this sort can be defined using the quaternionic, and
perhaps octonionic, points of a variety defined over $\Q$, but I have not pursued that either (yet).
\bigskip

\section{Genus zero moduli spaces} \bigskip

\noindent
{\bf 2.1} One way to explore the complex line at infinity is to send out a finite collection of
points as probes [18]: \bigskip

\noindent
The moduli spaces of stable genus zero curves, marked with suitably many distinct ordered
smooth points, have been extensively studied; they form a (slightly incomplete, in that $n \geq
2$) operad $\{ \bM_{0,n+1} \}$ of smooth algebraic varieties over $\Z$. Algebras over the
operad by the homology of its points (enlarged slightly to accommodate the Novikov ring)
are fundamental to the study of quantum cohomology: these are the polycommutative algebras of
Manin [8, 13]. On the other hand, the spaces $\{ \bM_{0,n+1}(\R) \}$ are acyclic, and their
fundamental groups define an operad in discrete groups [5, 6]. Recently, Etingof {\it et al} [7]
have identified algebras over $\{ H_*(\bM_{0,n+1}(\R),\Q) \}$ as 2-Gerstenhaber algebras: these
are commutative algebras endowed with a suitably compatible structure as a 2-Lie algebra (in the
sense of Hanlon and Wachs. Such a structure can be interpreted as a restricted kind of
$L_\infty$-algebra (carrying a triple bracket satisfying a generalized Jacobi identity), or [2 \S 2.3
Cor. 17] as a trivial Lie algebra together with an additional gerbey-looking 3-cohomology class.) 
One thus expects the existence of such a triple bracket on $H^*(V(\R),\Q[\pi_2V(\C)])$. \bigskip

\noindent
Note that the Hodge cohomology of an algebraic variety defined over the real is more than just
a Galois representation; it is naturally a representation of the Weil group of $\C/\R$ [21 \S 4.4].
\bigskip

\noindent
All these operads are in fact cyclic: the action of the symmetric group $\Sigma_n$ on
$\bM_{0,n+1}$ extends naturally to an action of $\Sigma_{n+1}$. My impression is that this
cyclic structure, on the real points, is still not completely understood. \bigskip

\noindent
{\bf 2.2} The real genus zero moduli spaces have a natural tesselation by Stasheff polyhedra, 
and thus a canonical CW-decomposition. There is a set $B_n$ of $(2n-3)!!$ distinguished zero-cells in
$\bM_{n+1}(\R)$, corresponding to configurations whose dual graphs are rooted binary planar
trees with $n$ labeled leaves, modulo a certain equivalence relation which allows twigs with
precisely two leaves to rotate [6, 19]. I am indebted to John Baez [1] for pointing out that 
this collection of zero - cells is itself an operad, whose algebras are commutative magmas 
(sets with a commutative but not necessarily associative binary product) and that its generating 
function
\[
|B|(T) = \sum |B_n| \frac{T^n}{n!} = 1 - (1 - 2T)^{1/2}
\]
satisfies an equation                                  
\[
|B|(T) = T + \half |B|(T)^2 
\]
which asserts that if an algebra over this operad is not a one-element set, then it is the union of an
unordered pair of sets, each with such a structure. \bigskip

\noindent
From another point of view such binary trees define `totally decomposed' configurations of genus
zero marked curves, with three points on each irreducible component; thus when $n=3$ we
recover configurations in $\bM_{0,4}$ corresponding to the three canonical points $0,1,\infty$
on $\PP_1$. In fact in general such maximally degenerate configurations are all defined over
$\Q$, and indeed  over $\Z$; they define a natural `operad of basepoints' for $\{\bM_{*+1}
(\Z)\}$. \bigskip

\noindent
{\bf 2.3} Applying the constructions of the section above defines a (cyclic) operad
\[
\{ \bPi \; (\bM_*(\C) \rel \bM_*(\R) \rel B_*) \}
\]
in two-categories.  Because $\bM_{0,n}(\R)$ has trivial higher homotopy groups, and because
$\bM_{0,n}(\C)$ has vanishing odd cohomology, there is an exact sequence
\[
0 \to \pi_2 (\bM_{0,n},\Z) \to \pi_2 (\bM_{0,n}(\C),\bM_{0,n}(\R)) \to \pi_1 (\bM_{0,n}(\R)) \to
1 
\]
(with trivial action of the quotient group, as far as I can see). This two-category then has Baez's
operad $B_*$ of commutative magmas (or, more precisely, its cyclic enrichment) as objects;
automorphisms of an object correspond to elements of the pure cactus group [identified
abstractly in [5] and in terms of generators and relations in [10 \S 3.4 ]], while the 
two-endomorphisms of a morphism correspond to elements of  $\pi_2 \cong H_2(\bM_{0,n}
(\C),\Z)$. According to Keel, this group is free abelian, with generators corresponding to
unordered partitions of $n$ into two sets, both of cardinality at least two, ie to trees with exactly
two internal nodes. \bigskip

\section{Coda: strings of pearls} \bigskip

\noindent
{\bf 3.1} A genus zero curve with only two marked points has too many automorphisms to be
stable, so the operad of DKM moduli spaces lacks a good `unary' operation. More generally,
Manin and Losev study stacks
\[
\BL_{g,m} = \coprod_{n \geq 1}  \bL_{g,m+1,n}
\]
of genus $g$ curves marked with $m+1$ marked smooth points, required to be distinct, and $n$
further auxiliary marked smooth points, which can collide with one another; they show that the
resulting moduli objects admit the clutching morphisms needed for the structure of a complete
operad. In particular, the monoid $\BL_{0,1}$ (which now has countably many components)
provides a replacement for the nonexistent $\bM_{0,2}$. \bigskip

\noindent
The configurations parametrized by this new space might be visualized as chains of bubbles,
exactly as seen in physicists' experiments, or as strings of pearls: sequences of $\PP_1$'s
laid end to end, with `heavy' marked points (corresponding to $0$ and $\infty$ ) anchoring the 
ends of the chain; intermediate nodes are marked by the ephemeral points of the second class. 
Just as the usual genus zero operad captures the WDVV equations for associativity, the resulting 
extended modular operad accounts for certain commutativity equations arising in physics [3]. 
\bigskip

\noindent
The operads $\BL_{0,*}$ and their cohomology have been extensively studied [15 \S 3.3.1, 16
\S 4.6.1] ; they, and their algebras, are relatively well-understood. In particular, the components
of $\bL_{0,2,*}$ are toric varieties, with permutohedra as associated quotients, and the
cohomology of $\bL_{0,2,n}$ has a combinatorical interpretation related to the Stanley - Reisner
ring of (the barycentric subdivision of) the $n$-cube [20 \S 4.3]. The real points $\BL_{0,1}(\R)$
are again aspherical [11] but as far as I know their fundamental groups have yet to be calculated
explicitly. \bigskip  

\noindent
{\bf 3.2} The construction of moduli stacks $\BL_{*,*}(V)$ of stable maps from such
configurations to some reasonable class of (algebraic, or more generally, symplectic) varieties, 
seems to be within reach [17]; another objective of this note is to draw attention to the interest of 
these objects. \bigskip

\noindent
$\BL_{0,1}(V)$, in particular, is a space over $V \times V$: it is an algebraic analog of the
classical space of (unrestriced) paths in $V$, and one expects the existence of glueing morphisms
making it into a monoid in such a category, with respect to the concatenation of bubble chains;
but (just as in the classical case) the relevant monoidal structure on spaces over $V \times V$ is
not symmetric. \bigskip

\noindent
If we fix a basepoint $\infty \in V$, then the fiber $\BL_{0,1}(V \rel \infty)$ of $\BL_{0,1}(V)$
over $(\infty,\infty) \in V \times V$ is the analog of the classical space of based loops in $V$,
which is an algebra over the $A_\infty$ operad. [The space of unbased loops is similarly an
algebra over this operad, in the category of spaces over $V$, with its natural symmetric monoidal
structure]. This construction is a good candidate for an algebra over $\{\BL_{0,*}\}$. If $V$ is
defined over $\Q$ or $\R$, this Manin - Losev loopspace will be as well, and its (rational)
homology can be expected to be a $W(\C/\R)$-equivariant $\{H_*(\BL_{0,*}\}$ - algebra.
This entails more than the existence of a nice action of the Galois group on the homology: in
particular it means that the various fixed point sets possess compatible operad actions. \bigskip

\noindent
Manin and Losev [15 \S 3.3.1] describe such algebras in an essentially Tannakian
language, so we may imagine thinking of them as representations of some motivic pro-algebraic
group. The remarks above suggest that it may be natural to enlarge this group slightly, to include
$Q \pi_1(\BL_{0,1})$ as something like its group of components.  \bigskip
           
\newpage

\bibliographystyle{amsplain}

\end{document}